\magnification=1200
\overfullrule=0pt
{\centerline {\bf Recent uses of connectedness in functional analysis}}\par
\bigskip
\bigskip
{\centerline {BIAGIO RICCERI}}\par
\bigskip
\bigskip
\bigskip
Perhaps, it is not too far from the truth to say that, among the great
concepts (as compactness, completeness, order, convexity) on which
functional analysis is based, connectedness is relatively less popular,
though this does not mean that it is less useful than the others. 
 The aim of this lecture is just to support this latter sentence,
focusing some recent results where connectedness plays a central role.\par
\smallskip
Our starting point is Theorem 1 below. Before stating it, to give
the reader the convenience to realize analogies and differences, we
recall, grouped together in Theorem A, three very famous results due
to K.Fan and F.E.Browder.\par
\smallskip
Given a product space
$X\times Y$, we denote by $p_{X}$ and $p_{Y}$ the projections from 
$X\times Y$ onto $X$ and $Y$, respectively. Moreover, if
$A\subseteq X\times Y$, for
every
 $x\in X$ and $y\in Y$, we put
$$A_{x}=\{v\in Y : (x,v)\in A\}$$
and
$$A^{y}=\{u\in X : (u,y)\in A\}.$$
\smallskip
THEOREM A ([5], Theorems 1 and 2; [1], Theorem 7). 
- {\it Let $E$, $F$ be two real Hausdorff 
locally convex topological vector spaces, let 
 $X\subseteq E$, $Y\subseteq F$ be two
 non-empty compact convex sets, and 
 let $S$, $T$ be two subsets
of $X\times Y$. Assume that at least one of the following three sets of
conditions is satisfied:\par
\noindent
$(\alpha)$\hskip 10pt  $S^{y}$ is convex
for each $y\in Y$, $S_{x}$ is open in Y for each $x\in X$, 
$T_{x}$ convex for each $x\in X$, 
and  $T^{y}$ is open in X for each $y\in Y$;\par
\noindent
$(\beta)$\hskip 10pt $S, T$ are closed, $S^{y}$ is convex
for each $y\in Y$, and $T_{x}$ is convex for each $x\in X$;
 \par
\noindent
$(\gamma)$\hskip 10pt  $S^{y}$ is convex
for each $y\in Y$, $S_{x}$ is open in Y for each $x\in X$, 
  $T$ is
closed, and $T_{x}$ is convex for each $x\in X$.\par
Then, at least one of the following assertions does hold:\par
\noindent
\hbox {\rm (a)}\hskip 10pt $p_{X}(T)\neq X.$\par
\noindent
\hbox {\rm (b)}\hskip 10pt $p_{Y}(S)\neq Y.$\par
\noindent
\hbox {\rm (c)}\hskip 10pt $S\cap T\neq\emptyset.$}\par
\medskip
In [18], we pointed out that, when $Y$ is a segment, Theorem A
is still true assuming simply that the sections $S^{y}$
are connected. More precisely, we have the following
\medskip
THEOREM 1 ([18], Theorem 2.3). - 
 {\it Let $X$, $Y$ be two topological spaces, with $Y$
admitting a continuous bijection onto $[0,1]$, and let $S$, $T$ be two
subsets of $X\times Y$, with  $S$ 
connected and, for each $x\in X$, $T_{x}$ connected.
Moreover, assume that either $T^{y}$ is open for each $y\in Y$, or
$Y$ is compact and $T$ is closed.\par
Then, at least one of the following assertions does hold:\par
\noindent
\hbox {\rm (a)}\hskip 10pt $p_{X}(T)\neq X.$\par
\noindent
\hbox {\rm (b)}\hskip 10pt $p_{Y}(S)\neq Y$ and $\{y\in Y :
 (p_{X}(S)\times
\{y\})\cap T = \emptyset\}\neq\emptyset.$\par
\noindent
\hbox {\rm (c)}\hskip 10pt $S\cap T\neq\emptyset.$}\par
\medskip
The following proposition is useful to recognize the connectedness
of a given set in a product space.\par
\medskip
PROPOSITION 1 ([18], Theorem 2.4). - {\it Let $X,Y$ be two topological
 spaces
 and let 
$S$ be a subset of $X\times Y$. Assume that at least one of the
following four sets of conditions is satisfied:\par
\noindent
$(\gamma_{1})$\hskip 10pt $p_{Y}(S)$ is connected, $S^{y}$ is connected
for each $y\in Y$, and $S_{x}$ is open for each $x\in X$;\par
\noindent
$(\gamma_{2})$\hskip 10pt $p_{Y}(S)$ is connected, $X$ is compact,
$S$ is closed, and $S^{y}$ is connected for each $y\in Y$;\par
\noindent
$(\gamma_{3})$\hskip 10pt $p_{X}(S)$ is connected, $S_{x}$ is connected
for each $x\in X$, and $S^{y}$ is open for each $y\in Y$;\par
\noindent
$(\gamma_{4})$\hskip 10pt $p_{X}(S)$ is connected, $Y$ is compact, $S$
is closed and $S_{x}$ is connected for each $x\in X$.\par
Under such hypotheses, $S$ is connected.}\par
\medskip
Then, thanks to Proposition 1, we have the following particular case of
Theorem 1 which is directly comparable with Theorem A (see also
[2]):\par
\medskip
THEOREM 2 ([18], Theorem 2.5). - {\it Let $X,Y$ be two topological spaces,
with $Y$ admitting a continuous bijection onto $[0,1]$, and let
$S,T$ be two subsets of $X\times Y$. Assume that at least one
of the following eight sets of conditions is satisfied:\par
\noindent
$(\delta_{1})$\hskip 10pt $p_{Y}(S)$ is connected, $S^{y}$ is connected
for each $y\in Y$, $S_{x}$ is open for each $x\in X$, 
$T_{x}$ is connected for each $x\in X$, 
and  $T^{y}$ is open for each $y\in Y$;\par
\noindent
$(\delta_{2})$\hskip 10pt $p_{Y}(S)$ is connected, $Y$ is
compact, $S^{y}$ is connected
for each $y\in Y$, $S_{x}$ is open for each $x\in X$, 
 $T$ is closed, and 
 $T_{x}$ is connected for each $x\in X$;\par
\noindent
$(\delta_{3})$\hskip 10pt $p_{Y}(S)$ is connected, $X$ is compact,
$S$ is closed, $S^{y}$ is connected for each $y\in Y$,
$T_{x}$ is connected for each $x\in X$, 
and  $T^{y}$ is open for each $y\in Y$;\par
\noindent
$(\delta_{4})$\hskip 10pt $p_{Y}(S)$ is connected, $X$ and $Y$ are compact,
$S$ and $T$ are closed, $S^{y}$ is connected for each $y\in Y$
, and 
 $T_{x}$ is connected for each $x\in X$;\par
\noindent
$(\delta_{5})$\hskip 10pt $p_{X}(S)$ is connected, $S_{x}$
and $T_{x}$ are connected for each $x\in X$,
and $S^{y}$ and $T^{y}$ 
 are open for each $y\in Y$; \par
\noindent
$(\delta_{6})$\hskip 10pt $p_{X}(S)$ is connected, $Y$ is compact,
 $S_{x}$ is connected for each $x\in X$,
$S^{y}$ is open for each $y\in Y$, $T$ is closed,
 and $T_{x}$ is connected for each $x\in X$;\par
 \noindent
$(\delta_{7})$\hskip 10pt $p_{X}(S)$ is connected, $Y$ is compact,
$S$ is closed, $S_{x}$ and $T_{x}$ are connected for each $x\in X$,
and $T^{y}$ is open for each $y\in Y$;\par
\noindent
$(\delta_{8})$\hskip 10pt $p_{X}(S)$ is connected, $Y$ is
compact, $S$ and $T$ are closed, and
 $S_{x}$ and $T_{x}$ are connected for each $x\in X$.\par
 Then, at least one of the following assertions does hold:\par
\noindent
\hbox {\rm (a)}\hskip 10pt $p_{X}(T)\neq X.$\par
\noindent
\hbox {\rm (b)}\hskip 10pt $p_{Y}(S)\neq Y$ and $\{y\in Y :
 (p_{X}(S)\times
\{y\})\cap T = \emptyset\}\neq\emptyset.$\par
\noindent
\hbox {\rm (c)}\hskip 10pt $S\cap T\neq\emptyset.$}\par
\medskip
We now start to present a series of applications of Theorems 1 and
2. The first of them concerns the following mini-max theorem:\par
\medskip
THEOREM 3 ([18], Theorem 1.1). -
  {\it Let $X,Y$ be two topological spaces, with $Y$ connected
and admitting a continuous bijection onto $[0,1]$, and
let $f$ be a real function on $X\times Y$. Assume that, for
each $\lambda>\sup_{y\in Y}\inf_{x\in X}f(x,y)$,
 $x_{0}\in X$, $y_{0}\in Y$, the sets
$$\{x\in X : f(x,y_{0})\leq \lambda\}$$
and
$$\{y\in Y : f(x_{0},y)>\lambda\}$$
are connected.
In addition, assume that at least one of the following three sets of
conditions is satisfied:\par
\noindent
$(h_{1})$\hskip 10pt $f(x,\cdot)$ is upper semicontinuous in
$Y$ for each $x\in X$, and $f(\cdot,y)$ is lower semicontinuous
in $X$ for each $y\in Y$;\par
\noindent
$(h_{2})$\hskip 10pt $Y$ is compact, and $f$ is upper semicontinuous
in $X\times Y$;\par
\noindent
$(h_{3})$\hskip 10pt $X$ is compact, and $f$ is lower semicontinuous
in $X\times Y$.\par
\noindent
Under such hypotheses, one has
$$\sup_{y\in Y}\inf_{x\in X}f(x,y)=\inf_{x\in X}\sup_{y\in Y}f(x,y).$$}
\smallskip
Two applications of Theorem 3 will be, in turn, presented later.\par
\smallskip
Another application of Theorem 2 yields the following result on the
existence of Nash equilibrium points which is directly comparable with
Theorem 4 of [5].\par
\medskip
THEOREM 4 ([21], Theorem 10). - {\it Let $X$ be a
 Hausdorff compact topological space,
 $Y$ an arc, and $f, g$ two continuous real functions
on $X\times Y$ such that, for each $\lambda\in {\bf R}$, $x_{0}\in
X$, $y_{0}\in Y$, the sets $\{x\in X : f(x,y_{0})\geq \lambda\}$
and $\{y\in Y : g(x_{0},y)\geq \lambda\}$ are connected.\par
Then, there exists $(x^{*},y^{*})\in X\times Y$ such that
$$f(x^{*},y^{*})=\max_{x\in X}f(x,y^{*})$$
and
$$g(x^{*},y^{*})=\max_{y\in Y}g(x^{*},y).$$}
\smallskip
Another consequence of Theorem 2 is the following\par
\medskip
THEOREM 5 ([21], Theorem 5). - {\it Let $E$ be an infinite-dimensional
 Hausdorff topological vector space $E$, $X\subseteq E$
a convex set with non-empty interior, $K\subseteq E$ a relatively
compact subset, $Y\subseteq {\bf R}$ a compact interval, and
$S,T$ two subsets of $X\times Y$. Assume that:\par
\noindent
\hbox {\rm (i)}\hskip 10pt
$S_{x}$ is open in $Y$ for each $x\in X\setminus K$, and
 $S^{y}$ is convex and with non-empty interior for each
$y\in Y$;\par
\noindent
\hbox {\rm (ii)}\hskip 10pt  $T_{x}$ is non-empty and connected for
 each $x\in X\setminus K$, and either $T^{y}\setminus K$
 is open in $X\setminus K$
 for each $y\in
Y$, or $Y$ is compact and $T\setminus (K\times Y)$ is closed
 in $(X\setminus K)\times Y$.\par
Then, for every set $V\subseteq X\times Y$ such
that $V^{y}$ is relatively compact in $E$ for each $y\in Y$
and $V_{x}$ is closed in $Y$ for each $x\in X\setminus K$, the
set $(S\setminus (V\cup (K\times Y)))\cap T$ is
non-empty.}
\medskip
Theorem 5 was applied in [3] by A.Chinn\`\i\  to obtain what seems
to be the first mini-max theorem, involving two functions $f, g$,
where it is not assumed that $f\leq g$. Her result is as follows:\par
\medskip
THEOREM 6 ([3], Theorem 1). - {\it  Let $E, X, K, Y$ be as in
Theorem 5, and let $f,g,h$ be three real functions on $X\times
Y$. Assume that:\par
\noindent
\hbox {\rm (a)}\hskip 10pt $f(x,\cdot)$ is quasi-concave in
$Y$ for each $x\in X\setminus K$, and either $f$ is upper
semicontinuous in $(X\setminus K)\times Y$ or $f(\cdot,y)$ is lower
semicontinuous in $X\setminus K$ for each $y\in Y$;\par
\noindent
\hbox {\rm (b)}\hskip 10pt $g(x,\cdot)$ is upper semicontinuous
in $Y$ for each $x\in X\setminus K$, and $g(\cdot,y)$ is upper
semicontinuous and quasi-convex in $X$ for each $y\in Y$;\par
\noindent
\hbox {\rm (c)}\hskip 10pt $h(x,\cdot)$ is upper semicontinuous
in $Y$ for each $x\in X\setminus K$, and the set
$\{x\in X : h(x,y)\geq \lambda\}$ is relatively compact in $E$
for each $y\in Y$ and each $\lambda> \sup_{v\in Y}\inf_{u\in X}
g(u,v)$;\par
\noindent
\hbox {\rm (d)}\hskip 10pt $f(x,y)\leq \max\{g(x,y),h(x,y)\}$
for each $(x,y)\in (X\setminus K)\times Y$. \par
Then, for every relatively compact set $H\subseteq E$, one
has
$$\inf_{x\in X\setminus H}\sup_{y\in Y}f(x,y)\leq
\sup_{y\in Y}\inf_{x\in X}g(x,y).$$}
\medskip
A joint application of Theorem 2 with the classical Mazurkiewicz
theorem on the covering dimension, yields Theorem 7 below which
could be of interest in control theory.\par
\smallskip
 Precisely,  
let $b$ be a positive real number and let $F$ be a given
multifunction from $[0,b]\times {\bf R}^{n}$ into ${\bf R}^{n}$.
We denote by ${\cal S}_{F}$ the set of all Carath\'eodory solutions
of the problem $x'\in F(t,x), x(0)=0$ in $[0,b]$. That is to
say
$${\cal S}_{F}=\{u\in AC([0,b],{\bf R}^{n}):
u'(t)\in F(t,u(t))\hskip 10pt \hbox {\rm a.e. in}\hskip 4pt
 [0,b],\hskip 3pt u(0)=0\}$$
where, of course, $AC([0,b],{\bf R}^{n})$ denotes the space
of all absolutely continuous functions from $[0,b]$ into
${\bf R}^{n}$. 
For each $t\in [0,b]$, put
$${\cal A}_{F}(t)=\{u(t) : u\in {\cal S}_{F}\}.$$
In other words, ${\cal A}_{F}(t)$ denotes the attainable
 set at time $t$. Also,
put
$$V_{F}=\bigcup_{t\in [0,b]}{\cal A}_{F}(t).$$
Finally, set
$$C_{F}=\{x\in {\bf R}^{n} : \{t\in [0,b] : x\in {\cal A}_{F}(t)\}
\hskip 7pt  \hbox {\rm is connected}\}.$$
\medskip
THEOREM 7 ([21], Theorem 9). - {\it Assume that $F$ 
has non-empty compact convex values and bounded 
range. Moreover, assume that $F(\cdot,x)$ is measurable
for each $x\in {\bf R}^{n}$ and that $F(t,\cdot)$ is
upper semicontinuous for a.e. $t\in [0,b]$.\par
 Then, for every non-empty 
connected set
$X\subseteq V_{F}\cap C_{F}$ which is open in its affine hull
and different from $\{0\}$, 
   one has the following alternative:\par
\noindent
 either
$$X\subseteq {\cal A}_{F}(b)$$ or
 $$\dim({\cal A}_{F}(t)\cap X)\geq \dim(X) - 1$$
for some $t\in ]0,b[$, where $\dim(X)$ denotes the covering dimension
of $X$.}\par
\medskip
It is also worth noticing another application of Theorem 2 which allowed
P.Cubiotti and B.Di Bella to get the following result, where
$\langle \cdot,\cdot \rangle$ denotes the usual inner product in
${\bf R}^{n}$.\par
\medskip
THEOREM 8 ([4], Theorem 4). - {\it Let $f:[0,1]\rightarrow
{\bf R}^{n}$ ($n\geq 2$) be a continuous function, and let
$Y=\{y\in {\bf R}^{n} : \parallel y\parallel =1\}$. Assume that,
for each $\sigma <0$, there exists $L_{\sigma}>0$ such that, for
each finite set $\{y_{1},...,y_{k}\}\subseteq Y$, there is a
set $\{t_{1},...,t_{k}\}\subseteq [0,1]$ such that
$$\langle f(t_{i}),y_{i}\rangle\geq \sigma\hskip 10pt \hbox
{\rm and}\hskip 10pt |t_{i}-t_{j}|\leq L_{\sigma}\parallel
y_{i}-y_{j}\parallel$$
for all $i,j=1,...,k$.\par
Then, $f$ vanishes at some point of $[0,1]$.}\par
\medskip
The next result comes out from a joint application of Theorem 1
with the classical Leray-Schauder continuation principle.\par
\medskip
THEOREM 9 ([21], Theorem 12). - {\it Let $E$ be a Banach space, $[a,b]$
a compact real interval, $\Omega$ a
non-empty open bounded subset of $E$, $f$ a continuous function
from $\overline {\Omega}\times [a,b]$ into $E$, with relatively
compact range. Assume that $f(x,y)\neq x$ for all
$(x,y)\in \partial\Omega\times [a,b]$ and that the
Leray-Schauder index of $f(\cdot,a)$ is not zero.\par
Then, for every lower semicontinuous function $\varphi:\Omega
\rightarrow [a,b]$ and every upper semicontinuous
function $\psi:\Omega\rightarrow [a,b]$, with $\varphi(x)\leq
\psi(x)$ for all $x\in \Omega$, there exist $x^{*}\in
\Omega$ and $y^{*}\in [\varphi(x^{*}),\psi(x^{*})]$ such that
$f(x^{*},y^{*})=x^{*}.$\par
In addition, if for some sequence $\{\lambda_{n}\}$ of positive real
numbers,
 with $\inf_{n\in {\bf N}}\lambda_{n}=0$, one has
$$\inf\{y\in [a,b] : \parallel f(x,y)-x\parallel\geq \lambda_{n}\}=
\inf\{y\in [a,b] : \parallel f(x,y)-x\parallel> \lambda_{n}\}$$
for each $x\in \Omega$, $n\in {\bf N}$ for
which
$$\{y\in [a,b] : \parallel f(x,y)-x\parallel >\lambda_{n}\}\neq
\emptyset,$$
 then there exists $x_{0}\in
\Omega$ such that $f(x_{0},y)=x_{0}$ for all $y\in [a,b]$.}\par
\medskip
We now come to the two announced applications of Theorem 3. 
The first of them is due to O.Naselli ([8]). Making use of Theorem
3, she got, as a corollary of a more general result, the following\par
\medskip
THEOREM 10 ([8], Theorem 3.4). - {\it Let $E$ be a
real Hausdorff topological vector
space, $p$ a real number greater than $1$, and $\alpha, \beta,
\gamma$ three affine functionals on $E$, with $\gamma(0)\geq 0$.\par
Then, for every closed, bounded and convex set $X\subseteq
\gamma^{-1}([\gamma(0),+\infty[)\cap \gamma^{-1}
(]0,+\infty[)$, with $\dim(X)\geq 2$,
one has
$$\inf_{x\in X}\left (\alpha(x)+\left ({{|\beta(x)|^{p}}
\over {\gamma(x)}}\right )^{1\over p-1}\right )=
\inf_{x\in B_{X}}\left (\alpha(x)+\left ({{|\beta(x)|^{p}}
\over {\gamma(x)}}\right )^{1\over p-1}\right ),$$
where
$$B_{X}=\{x\in X : \exists y\in \hbox {\rm aff}(X)\setminus 
\{x\} : [x,y]\cap X=\{x\}\},$$
$\hbox {\rm aff}(X)$ being the affine hull of $X$, and
$[x,y]$ being the line segment joining $x$ and $y$.}\par
\medskip
The other application of Theorem 3 we wish to recall concerns
integral functionals. We first introduce some notation.\par
\smallskip
 In the next four results,
  $(T,{\cal F},\mu)$ is a $\sigma$-finite
non-atomic measure space ($\mu(T)>0$), 
 $(E,\parallel
\cdot\parallel)$ is a real Banach space ($E\neq \{0\}$),
 and $p$ is a
real number greater than or equal to $1$. When $p=1$, we
will adopt the convention ${{p}\over {p-1}}=\infty$.\par
\smallskip
 For simplicity,
we denote by $X$ the usual space
 $L^{p}(T,E)$  of (equivalence
classes of) strongly $\mu$-measurable functions $u : T\rightarrow E$ 
such that
$\int_{T}\parallel u(t)\parallel^{p} d\mu<+\infty$, equipped with
the norm $\parallel u\parallel_{X}=
(\int_{T}\parallel u(t)\parallel^{p}d\mu)^{1\over p}.$\par
\smallskip
Moreover, we denote by ${\cal V}(X)$ the family of all sets
$V\subseteq X$ of the following type:
$$V=\{u\in X : \Psi(u)=\int_{T}g(t,u(t))d\mu\}$$
where $\Psi$ is a continuous linear
functional on $X$, and $g : T\times E\rightarrow {\bf R}$ is
such that the integral functional $u\rightarrow \int_{T}g(t,u(t))d\mu$
is (well-defined and) Lipschitzian in $X$, with Lipschitz constant
strictly less than $\parallel \Psi\parallel_{X^{*}}$.\par
\smallskip
Note, in particular, that each closed hyperplane of $X$ belongs to the 
family ${\cal V}(X)$.\par
\smallskip
 We then have\par
\medskip 
THEOREM 11 ([22], Theorem 2). - {\it Let $f:T\times
 E\rightarrow [0,+\infty[$ be
such that $f(\cdot,x)$ is $\mu$-measurable for each $x\in E$ and
$f(t,\cdot)$ is Lipschitzian
 with Lipschitz constant $M(t)$ for almost every $t\in T$, where
$M\in L^{p\over p-1}(T)$.
 Assume that $f(\cdot,0)\in L^{1}(T)$ and that
there exists
 a sequence $\{\lambda_{n}\}$ in $]0,+\infty[$,
with $\lim_{n\rightarrow +\infty}\lambda_{n}=+\infty$,
 such that,
 for almost
every $t\in T$ and for every $x\in E$, one has
$$\lim_{n\rightarrow +\infty}{{f(t,\lambda_{n} x)}\over
 {\lambda_{n}}}=0.$$
Then, for every $V\in {\cal V}(X)$, one has
$$\inf_{u\in V}\int_{T}f(t,u(t))d\mu=\inf_{u\in X}
\int_{T}f(t,u(t))d\mu.$$}\par
\smallskip
The proof of Theorem 11 is fully based on an application of
Lemma 1 of [19]. It is just this latter to be obtained by means
of an application of Theorem 3. It is also worth noticing that
such an application is made possible by the following very interesting
result of J.Saint Raymond:\par
\medskip
THEOREM 12 ([23], Th\'eor\`eme 3). - {\it Let $f:T\times E\rightarrow
{\bf R}$ be a ${\cal F}\otimes {\cal B}(E)$-measurable function,
${\cal B}(E)$ being the Borel family of $E$. Then, if we put
$$Y=\{u\in X : f(\cdot,u(\cdot))\in L^{1}(T)\},$$
for each $\lambda\in {\bf R}$, the set
$$\{u\in Y : \int_{T}f(t,u(t))d\mu\leq \lambda\}$$
is connected.}\par
\medskip
Theorem 11 has the following two consequences.\par
\medskip
THEOREM 13 ([22], Theorem 1). - {\it Let $E$ be separable, and let
 $F:T\rightarrow 2^{E}$ be a measurable multifunction,
with non-empty closed values. Assume that
 $\hbox {\rm dist}(0,F(\cdot))\in L^{1}(T)$
and that there exists a sequence $\{\lambda_{n}\}$ in $]0,+\infty[$,
with $\lim_{n\rightarrow +\infty}\lambda_{n}=+\infty$,
 such that,
 for almost every $t\in T$ and for every $x\in E$, one has
$$\lim_{n\rightarrow +\infty}
{{\hbox {\rm dist}(\lambda_{n} x,F(t))}\over {\lambda_{n}}}=0.$$
Then, if $p=1$, each member of the family ${\cal V}(X)$ contains
a selection of $F$.}\par
\medskip
 THEOREM 14 ([22], Theorem 6). - {\it Let $E$ be reflexive and 
separable, let $p>1$,
 and let $f:T\times E\rightarrow [0,+\infty[$ be such that
$f(\cdot,x)$ is $\mu$-measurable for each $x\in E$,
$f(\cdot,0)\in L^{1}(T)$,
 and $f(t,\cdot)$ is G\^{a}teaux differentiable 
for almost every $t\in T$. Moreover, assume that there exist
$M\in L^{p\over p-1}(T)$ and
 a sequence
$\{\lambda_{n}\}$ in $]0,+\infty[$, with 
$\lim_{n\rightarrow +\infty}\lambda_{n}=+\infty$,
 such that, for almost every $t\in T$ and for
every $x\in E$, one has 
$$
\parallel
f_{x}'(t,x)\parallel_{E^{*}}\leq M(t)$$ and
$$\lim_{n\rightarrow +\infty}
{{f(t,\lambda_{n} x)}\over {\lambda_{n}}}=0.$$
Then, for every $V\in {\cal V}(X)$, there exists a
sequence $\{u_{n}\}$ in $V$ such that 
$$\lim_{n\rightarrow +\infty}\int_{T}f(t,u_{n}(t))d\mu=\inf_{u\in X}
\int_{T}f(t,u(t))d\mu$$
and
$$\lim_{n\rightarrow +\infty}\int_{T}\parallel 
f_{x}'(t,u_{n}(t))\parallel_{E^{*}}^{p\over p-1}d\mu=0.$$}\par
\smallskip
The final part of our lecture is devoted to recent applications
of the following lower semicontinuity result, based itself on connectedness:
\par
\medskip
THEOREM 15 ([10], Th\'eor\`eme 1.1). - {\it Let $X, Y$ be two topological
 spaces, with $Y$
connected and locally connected, and let $\varphi : X\times Y
\rightarrow {\bf R}$ be a function satisfying the following two
conditions:\par
\noindent
\hbox {\rm (a)}\hskip 10pt for each $x\in X$, the function
$\varphi(x,\cdot)$ is continuous, $0\in \hbox {\rm int}(\varphi(x,Y))$,
and $\hbox {\rm int}(\{y\in Y : \varphi(x,y)=0\})=\emptyset$;\par
\noindent
\hbox {\rm (b)}\hskip 10pt the set
$$\{(y,z)\in Y\times Y : \{x\in X : \varphi(x,y)<0<\varphi(x,z)\}
\hskip 5pt is\hskip 5pt open\}$$
is dense in $Y\times Y$.\par
Then, if, for each $x\in X$, one denotes by $Q(x)$ the set of
all $y\in Y$ such that $\varphi(x,y)=0$ and $y$ is not a local
extremum for $\varphi(x,\cdot)$, one has that $Q(x)$ is non-emtpy
and closed, and that the multifunction $x\rightarrow Q(x)$ is
lower semicontinuous.}
\medskip
We now recall two applications of Theorem 15 to implicit differential
equations.\par
\medskip
THEOREM 16 ([17], Theorem 2).
 - {\it Let $Y$ be a linear subspace of ${\bf R}^{n}$,
with $\dim(Y)\geq 2$, and let $f:[0,1]\times {\bf R}^{nk}\times
Y\rightarrow
{\bf R}$ be a continuous functions such that, for each 
$(t,\xi)\in [0,1]\times {\bf R}^{nk}$, $f(t,\xi,\cdot)$ is
affine and non-constant in $Y$.\par
Then, for every $x_{0}, x_{1},...,x_{k-1}\in {\bf R}^{n}$, there
exists $b\in ]0,1]$ such that the set of all functions
$u\in C^{k}([0,b],{\bf R}^{n})$ satisfying
$$u^{(k)}(t)\in Y, \hskip 3pt f(t,u(t),u'(t),...,u^{(k)}(t))=0
\hskip 5pt in\hskip 5pt [0,b],$$ 
$$u^{(i)}(0)=x_{i}\hskip 5pt for\hskip 5pt i=0,1,...,k-1,$$
has the continuum power.}\par
\medskip
THEOREM 17 ([7], Example 4.1). - {\it Let $\Omega\subseteq
{\bf R}^{n}$ ($n\geq 3$) be an open, bounded, connected
subset, with a boundary of class $C^{1,1}$.\par
Then, for every $g\in L^{p}(\Omega)$, with $p\in ]n,+\infty[$ , $\gamma
\in [0,1[$, $\lambda, \mu\in {\bf R}$, there exists
$u\in W^{2,p}(\Omega)\cap W^{1,p}_{0}(\Omega)$ such
that
$$\Delta u(x)=\lambda \sin \Delta u(x)+\mu
(|u(x)|+\parallel \nabla u(x)\parallel)^{\gamma}+
g(x)$$
for almost every $x\in \Omega$.}\par
\medskip
For other papers related to Theorem 15, we refer to [9],
[11], [12], [13], [14], [15].\par
\smallskip
Before establishing the final applications of Theorem 15, we also
recall the following\par
\medskip
THEOREM 18 ([16], Th\'eor\`eme 2). 
 - {\it Let $X, Y$ be two real Banach spaces, let
$\Phi:X\rightarrow Y$ be a continuous linear surjective operator,
and let $\Psi:X\rightarrow Y$ be a Lipschitzian operator, with
Lipschitz constant $L<{{1}\over {\alpha_{\Phi}}}$, where
$\alpha_{\Phi}=\sup_{\parallel y\parallel\leq 1}\hbox
 {\rm dist}(0,\Phi^{-1}(y)).$\par
Then, for each $y\in Y$, the set $(\Phi+\Psi)^{-1}(y)$ is a
(non-empty) retract of $X$, and the multifunction
$y\rightarrow (\Phi+\Psi)^{-1}(y)$ is Lipschitzian (with respect
to the Hausdorff distance), with Lipschitz constant 
${{\alpha_{\Phi}}\over {1-L\alpha_{\Phi}}}.$}\par
\medskip
We now can prove\par
\medskip
THEOREM 19. - {\it Let $X$ be a connected topological space,
 $E$ a real Banach space (with topological dual space $E^{*}$),
  $\Phi$ an operator from $X$ into $E^{*}$, $f$
a real function on $X\times E$ such that,
for each $x\in X$, $f(x,\cdot)$ is Lipschitzian in $E$, with 
Lipschitz constant $L(x)\geq 0$. Further, assume that the
set 
$$\{y\in E : \langle \Phi(\cdot),y\rangle - f(\cdot,y)\hskip 5pt
is\hskip 5pt continuous\}$$ 
is dense in $E$ and that the set
$$\{(x,y)\in X\times E : \langle \Phi(x),y\rangle =
f(x,y)\}$$
is disconnected.\par
Then, there exists some $x_{0}\in X$ such that
$\parallel \Phi(x_{0})\parallel_{E^{*}}\leq L(x_{0}).$}\par
\smallskip
PROOF. Arguing by contradiction, assume that 
$\parallel \Phi(x)\parallel_{E^{*}}> L(x)$ for all $x\in X$.
 Then, by Theorem 18, for each $x\in X$, the function $\langle
\Phi(x),\cdot\rangle - f(x,\cdot)$ is onto ${\bf R}$, is open and
has connected point inverses. At this point, we can apply    
Theorem 15, to get that
the multifunction $Q:X\rightarrow 2^{E}$ defined
by
$$Q(x)=\{y\in E : \langle \Phi(x),y\rangle = f(x,y)\}$$
 is lower semicontinuous. Then, since
$X$ is connected and each $Q(x)$ is non-empty and connected, Theorem
3.2 of [6] ensures that the graph of $Q$ is connected too, against
one of our assumptions.\hskip 10pt $\bigtriangleup$
\par
\medskip
Observe that when, in Theorem 19, $f$ does not depend on $y$ (that is,
$L(x)=0$ for all $x\in X$)
we directly get the existence of a zero for the operator $\Phi$. In this
case, one can even assume that $E$ is simply a topological vector
space (see [20]). To get a zero for $\Phi$ allowing $f$ to depend on
$y$, we can use\par
\medskip
THEOREM 20. - {\it Let $X$ be a connected topological space,
 $E$ a real Banach space (with topological dual space $E^{*}$),
  $\Phi$ an operator from $X$ into $E^{*}$, with closed range. Assume
that, for each $\epsilon>0$, there exists a real function $f_{\epsilon}$ on
$X\times E$ having the following properties:\par
\noindent
\hbox {\rm (a)}\hskip 10pt for each $x\in X$, the function
$f_{\epsilon}(x,\cdot)$ is Lipschitzian in $E$, with Lipschitz constant
 less than or equal to $\epsilon$;\par
\noindent
\hbox {\rm (b)}\hskip 10pt the set 
$$\{y\in E : \langle \Phi(\cdot),y\rangle - f_{\epsilon}(\cdot,y)
\hskip 5pt is\hskip 5pt continuous\}$$ 
is dense in $E$;\par
\noindent
\hbox {\rm (c)}\hskip 10pt   the set
$$\{(x,y)\in X\times E : \langle \Phi(x),y\rangle =
f_{\epsilon}(x,y)\}$$
is disconnected.\par
Then, $\Phi$ vanishes at some point of $X$.}\par
\smallskip
PROOF. Applying Theorem 19, for each $\epsilon>0$, we get a point
$x_{\epsilon}\in X$ such that $\parallel \Phi(x_{\epsilon})
\parallel_{X^{*}}\leq \epsilon$. In other words, $0$ is in the
closure of $\Phi(X)$. But, by assumption, $\Phi(X)$ is closed,
and so $0\in \Phi(X)$, as claimed.\hskip 10pt $\bigtriangleup$\par
\medskip
THEOREM 21. - {\it Let $X$ be a connected and locally connected
topological space, $E$ a real Banach space, $\Phi : X\rightarrow
E^{*}$ a (strongly) continuous operator, $L$ a non-negative real
function on
$X$. Denote by $\Lambda$ the set of all continuous functions $f:X\times E
\rightarrow {\bf R}$ such that, for each $x\in X$, $f(x,\cdot)$
is Lipschitzian in $E$, with Lipschitz constant less than or equal
to $L(x)$.
Consider $\Lambda$ equipped with the relativization of the
strongest vector topology on the space ${\bf R}^{X\times E}$, and
assume that the set
$$\{(f,x,y)\in \Lambda\times X\times E : \langle \Phi(x),y\rangle =
f(x,y)\}$$
is disconnected.\par
Then, there exists some $x_{0}\in X$ such that
$\parallel \Phi(x_{0})\parallel_{E^{*}}\leq L(x_{0}).$}\par
\smallskip
PROOF. Arguing by contradiction, assume that 
$\parallel \Phi(x)\parallel_{E^{*}}> L(x)$ for all $x\in X$. For
each $(f,x,y)\in \Lambda\times X\times E$, put
$$\varphi(f,x,y)=\langle \Phi(x),y\rangle - f(x,y).$$
Observe that, for each $(x,y)\in X\times E$, the function
$\varphi(\cdot,x,y)$ is continuous in $\Lambda$ since it is 
continuous even with respect to the topology of pointwise convergence.
Moreover, for each $f\in \Lambda$, the function $\varphi(f,\cdot,
\cdot)$ is continuous in $X\times E$ (thanks to the strong continuity
of $\Phi$), and, again by Theorem 18, is onto ${\bf R}$, and
has no local extrema. Consequently, again by Theorem 15,
the multifunction $Q:\Lambda\rightarrow 2^{X\times E}$ defined by
$$Q(f)=\{(x,y)\in X\times E : \langle \Phi(x),y\rangle =
f(x,y)\}$$
is lower semicontinuous. But, by Theorem 19, each set $Q(f)$ is
connected. On the other hand, $\Lambda$ is connected (since it is
convex), and so, by Theorem 3.2 of [6], the graph of $Q$ is connected,    
against one of our assumptions.\hskip 10pt $\bigtriangleup$
\medskip
>From Theorem 21, we get, of course, the following\par
\medskip
THEOREM 22. - {\it Let $X$ be a connected and locally connected
topological space, $E$ a real Banach space, $\Phi : X\rightarrow
E^{*}$ a (strongly) continuous operator, with closed range. For
each $\epsilon>0$, denote by $\Lambda_{\epsilon}$
the set of all continuous functions $f:X\times E
\rightarrow {\bf R}$ such that, for each $x\in X$, $f(x,\cdot)$
is Lipschitzian in $E$, with Lipschitz constant less than or equal
to $\epsilon$.
Consider $\Lambda_{\epsilon}$ equipped with the relativization of the
strongest vector topology on the space ${\bf R}^{X\times E}$, and
assume that the set
$$\{(f,x,y)\in \Lambda_{\epsilon}\times X\times E :
 \langle \Phi(x),y\rangle =f(x,y)\}$$
is disconnected.\par
Then, $\Phi$ vanishes at some point of $X$.}\par
\medskip
We conclude by proposing the following conjecture:\par 
\medskip
CONJECTURE. {\it Let $X$ be the closed unit ball of ${\bf R}^{n}$,
$(n\geq 2)$, $g:X\rightarrow X$ a continuous function and $\epsilon>0$.
Denote by $\Lambda_{\epsilon}$ the set all continuous
functions $f:X\times {\bf R}^{n}
\rightarrow {\bf R}$ such that, for each $x\in X$,
 $f(x,\cdot)$ is Lipschitzian in
${\bf R}^{n}$, with Lipschitz constant less than or equal to
$\epsilon$. Consider $\Lambda_{\epsilon}$ equipped with the relativization of
the strongest vector topology on the space ${\bf R}^{X\times {\bf R}^{n}}$.
  Then, the set
$$\{(f,x,y)\in \Lambda_{\epsilon}
\times X\times {\bf R}^{n} : \langle g(x)-x,y\rangle
=f(x,y)\}$$
is disconnected.}\par
\medskip
Observe that, on the basis of Theorem 22, the above conjecture could
lead to a completely new way of proving the Brouwer fixed point theorem.
\par
\bigskip
\bigskip
{\centerline {\bf References}}
\bigskip
\bigskip
\noindent
[1]\hskip 5pt F.E.BROWDER, {\it The fixed point theory of multi-valued
mappings in topological vector spaces}, Math.Ann., {\bf 177} (1968),
283-301.\par
\smallskip
\noindent
[2]\hskip 5pt A.CHINN\`I, {\it Some remarks on a theorem on sets with
connected sections}, Rend.Circ.Mat.\par
\noindent
Palermo, to appear.\par
\smallskip
\noindent
[3]\hskip 5pt A.CHINN\`I, {\it A two-function minimax theorem},
Arch.Math.(Basel), to appear.\par
\smallskip
\noindent
[4]\hskip 5pt P.CUBIOTTI and B.DI BELLA, {\it A $\sup-\inf$-condition
for the existence of zeros of certain nonlinear operators}, 
Ann.Univ.Sci.Budapest E\"otv\"os Sect.Math., to appear.\par
\smallskip
\noindent
[5]\hskip 5pt K.FAN, {\it Applications of a theorem concerning sets
with convex sections}, Math.Ann., {\bf 163} (1966), 189-203.\par
\smallskip
\noindent
[6]\hskip 5pt J.-B.HIRIART-URRUTY, {\it Images of connected sets by
semicontinuous multifunctions}, J.Math.Anal.Appl., {\bf 111} (1985),
407-422.\par
\smallskip
\noindent
[7]\hskip 5pt S.A.MARANO, {\it Implicit elliptic differential
equations}, Set-Valued Anal., {\bf 2} (1994), 545-558.\par
\smallskip
\noindent
[8]\hskip 5pt O.NASELLI, {\it On a class of functions with
equal infima over a domain and over its boundary}, J.Optim.Theory
Appl., to appear.\par
\smallskip
\noindent
[9]\hskip 5pt B.RICCERI, {\it Sur la semi-continuit\'e inf\'erieure
de certaines multifonctions}, C.R.Acad.Sci.\par
\noindent
Paris, S\'erie I,
{\bf 294} (1982), 265-267.\par
\smallskip
\noindent
[10]\hskip 5pt B.RICCERI, {\it Applications de th\'eor\`emes de
semi-continuit\'e inf\'erieure}, C.R.Acad.Sci.Paris, S\'erie I,
{\bf 295} (1982), 75-78.\par
\smallskip
\noindent
[11]\hskip 5pt B.RICCERI, {\it Solutions lipschitziennes
d'\'equations diff\'erentielles sous forme implicite}, C.R.Acad.
Sci.Paris, S\'erie I, {\bf 295} (1982), 245-248.\par
\smallskip
\noindent
[12]\hskip 5pt B.RICCERI, {\it On multiselections}, Matematiche,
{\bf 38} (1983), 221-235.\par
\smallskip
\noindent
[13]\hskip 5pt B.RICCERI, {\it On inductively open real functions},
Proc.Amer.Math.Soc., {\bf 90} (1984), 485-487.\par
\smallskip
\noindent
[14]\hskip 5pt B.RICCERI, {\it Lifting theorems for real functions},
Math.Z., {\bf 186} (1984), 299-307.\par
\smallskip
\noindent
[15]\hskip 5pt B.RICCERI, {\it On multifunctions of one real variable},
J.Math.Anal.Appl., {\bf 124} (1987), 225-236.\par
\smallskip
\noindent
[16]\hskip 5pt B.RICCERI, {\it Structure, approximation et d\'ependance
continue des solutions de certaines \'equations non lin\'eaires},
C.R.Acad.Sci.Paris, S\'erie I, {\bf 305} (1987), 45-47.\par
\smallskip
\noindent
[17]\hskip 5pt B.RICCERI, {\it On the Cauchy problem for the differential
equation $f(t,x,x',...,x^{(k)})=0$}, Glasgow Math.J., {\bf 33} (1991),
343-348.\par
\smallskip
\noindent
[18]\hskip 5pt B.RICCERI, {\it Some topological mini-max theorems via
an alternative principle for multifunctions}, Arch.Math.(Basel),
{\bf 60} (1993), 367-377.\par
\smallskip
\noindent
[19]\hskip 5pt B.RICCERI, {\it A variational property of integral 
functionals on $L^{p}$-spaces of vector-valued functions}, 
C.R.Acad.Sci.Paris, S\'erie I, {\bf 318} (1994), 337-342.\par
\smallskip
\noindent
[20]\hskip 5pt B.RICCERI, {\it Existence of zeros via disconnectedness},
J.Convex Anal., {\bf 2} (1995), 287-290.\par
\smallskip
\noindent
[21]\hskip 5pt B.RICCERI, {\it Applications of a theorem concerning sets
with connected sections}, Topol.\par \noindent Methods Nonlinear Anal.,
to appear.\par
\smallskip
\noindent
[22]\hskip 5pt B.RICCERI, {\it On the integrable selections of certain
multifunctions}, Set-Valued Anal., to appear.\par
\smallskip
\noindent
[23]\hskip 5pt J.SAINT RAYMOND, {\it Connexit\'e des sous-niveaux
des fonctionnelles int\'egrales}, Rend.\par
\noindent
Circ.Mat.Palermo, {\bf 44} (1995), 162-168.\par
\bigskip
\bigskip
\noindent
Department of Mathematics\par
\noindent
University of Catania\par
\noindent
Viale A.Doria 6\par
\noindent
95125 Catania, Italy
\bye